\input amstex
\documentstyle {amsppt}
\input epsf
\advance\voffset  by -1.0cm
\NoBlackBoxes
\magnification=\magstep1
\hsize=17truecm
\vsize=24.2truecm
\voffset=0.5truecm
\document

\define \Ga {\Gamma}
\define \la {\lambda}
\define \pa {\partial}
\define \bR {\Bbb R}
\define \bC {\Bbb C}
\define \al {\alpha}

\def \sp {\operatorname{supp}}

\define \polr {P_{\bR}}
\define \polc {P_{\bC}}

\topmatter
\title
On Bochner-Krall orthogonal polynomial systems
\endtitle

\author T.~Bergkvist, H.~Rullg\aa rd and B.~Shapiro
\endauthor
\affil
   Department of Mathematics, University of Stockholm, S-10691, Sweden,
   {\tt tanjab\@matematik.su.se},
{\tt hansr\@matematik.su.se},  {\tt shapiro\@matematik.su.se}\\
\endaffil
\rightheadtext {On Bochner-Krall OPS}
\leftheadtext {Bergkvist, Rullg\aa rd and Shapiro}

\abstract
In this paper we address the classical question going back to
S.~Bochner and H.~L.~Krall (see original \cite {B}, \cite{Kr} and \cite {EKLW}
for the modern status quo) to describe
all systems  $\{p_{n}(x)\}_{n=0}^\infty$ of orthogonal polynomials (OPS)
which are the eigenfunctions  of some finite order differential
operator, i.e.  satisfy  the equation
$$\sum_{k=1}^{N}a_{k}(x)y^{(k)}(x)=\la_{n} y(x).\tag1$$
Such systems of orthogonal polynomials are  called {\it Bochner-Krall OPS}
(or BKS for short) and their spectral differential operators are
accordingly called
{\it Bochner-Krall operators} (or BK-operators for short).
We say that a BKS has compact type if it is orthogonal with respect to a
compactly supported positive measure on $\bR$.
   It is well-known that the order $N$ of any BK-operator should be even
   and every coefficient $a_{k}(x)$ must be a
polynomial of degree at most $k$, see e.g. \cite {EKLW}.
   Below we show that the leading coefficient of a compact type 
   BK-operator is of the form $((x - a)(x-b))^{N/2}$.
   This settles the special case of  the
general conjecture 7.3. of
   \cite {EKLW}
   describing the leading terms of all BK-operators. New results on the
   asymptotic distribution of zeros of polynomial eigenfunctions
   for a spectral problem (1)  are the main ingredient
   in the proofs, comp. \cite {BR} and \cite
   {MS}.

\endabstract

\endtopmatter

\heading \S 1. Introduction  \endheading
Let $\polr$ and $\polc$ denote the spaces of all real and complex
polynomials, respectively, in a variable
$x$. By a real (or complex) {\it polynomial system} we will mean
a sequence $\{p_{n}\}_{n = 0}^\infty$ of polynomials in
$\polr$ (or $\polc$) such that $\deg p_{n} = n$. By an orthogonal
polynomial system (OPS) one understands a real polynomial system
$\{p_{n}\}$ such that $\langle p_{n}, p_{m} \rangle$ is nonzero
precisely when $n = m$, where $\langle \, , \, \rangle$ is some
reasonable inner product on the linear space $\polr$.
If an orthogonal polynomial system for a given inner product exists, the
$p_{n}$ are unique, up to multiplication by scalars.
\medskip
Orthogonal polynomial systems have been studied in various degrees of
generality. Classically, one has considered inner products of the form
$$
\langle p, q \rangle = \sigma(p\cdot q)
$$
where $\sigma$ is a moment functional, that is a linear functional on $\polr$.
It is known that all moment functionals can be represented by an integral
$$
     \sigma(p) = \int p(x)\,d\mu(x)
$$
where $\mu$ is a (possibly signed) Borel measure on the real line.
The most complete theories have been obtained in the case
where $\mu$ is positive and has compact support.
In what follows we will mainly be concerned with orthogonal polynomial
systems of this particular kind, which we will call {\it 
compact type orthogonal systems }. Recently, there has been interest in 
more general inner
products called Sobolev, which are of the form
$$
\langle p, q \rangle = \sum_{k = 0}^{M} \sigma_{k}(p^{(k)} \cdot
q^{(k)})
$$
where the $\sigma_{k}$ are moment functionals. For the basics of the 
classical theory of
orthogonal polynomials see e.g. \cite {Sz} and \cite {Ch}.
\medskip
Consider now a differential operator
$$
     \frak d=\sum_{k=1}^{N}a_{k}(x)\frac {d^k}{dx^k} \tag 2
$$
where the
coefficients $a_{k}(x)$ are polynomials in $\polc$. We are interested in
eigenpolynomials of this operator, that is polynomials $p \in \polc$
satisfying $\frak d p = \lambda p$ for some constant $\lambda$.
Already S.~Bochner observed  that
the operator $\frak d$ has infinitely many
linearly independent eigenpolynomials if and only if $\deg a_{k} \leq
k$, with equality for at least one $k$. In this case there is
  precisely one monic degree $n$ eigenpolynomial
$p_{n}$ for all sufficiently large $n$. For generic $a_{k}$, the same
is true for every $n \geq 0$. If a complex polynomial system consists of
eigenpolynomials for an operator of the form (2), we will call it a
{\it system of eigenpolynomials}.
\medskip
A Bochner-Krall system (BKS for short) is defined to be a real polynomial
system which is both orthogonal (with respect to some inner product
$\langle \, , \, \rangle$)
and a system of eigenpolynomials (for some differential operator
$\frak d$). In this case $\langle \, , \, \rangle$
is called a Bochner-Krall inner
product, and $\frak d$ is called a Bochner-Krall operator.
If a BKS is a compact type orthogonal system, we will call it for short a
{\it compact type BKS}. The results we report in this note are valid for 
all compact type 
BKS.
\medskip
It is an open problem to classify all Bochner-Krall systems. A
complete classification is only known for Bochner-Krall operators
$\frak d$ with $N \leq 4$. The corresponding BKS are various classical
systems such as the Jacobi-type, the Laguerre-type, the Legendre 
type, the Bessel and the Hermite polynomials, see
\cite{EKLW}, Th. 3.1. In general, it is not even known which differential
operators are Bochner-Krall operators for some BKS.
In \cite{EKLW} it is conjectured that the leading
coefficient $a_{N}$ of a Bochner-Krall operator is a power of a
polynomial of degree at most 2. Our main result is an affirmative
answer to this conjecture for compact type BKS.
\medskip
Our results are obtained by studying the asymptotic distribution of
zeros of a polynomial system. To make
this more precise, let $\{p_{n}\}$ be a polynomial system
and for fixed $n \geq 1$, let $\alpha_{1}, \ldots, \alpha_{n}$
denote the (real or complex) zeros of $p_{n}$. Let $\nu_{n} =
\frac{1}{n} \sum_{i = 1}^n \delta(x - \alpha_{i})$ be the probability
measure in the complex plane with point masses at these zeros. We call
the measures $\nu_{n}$ the {\it root measures} of the polynomial system
$\{p_{n}\}$. If the sequence of root measures $\nu_{n}$ converges
weakly to a measure $\nu$ when $n \to \infty$, we say that $\nu$ is
the asymptotic distribution of zeros of the polynomial system.
\medskip
The following results, which characterize the asymptotic distribution
of zeros for compact type OPS and for systems of
eigenpolynomials respectively, are crucial to our treatment.
\medskip
Suppose that a polynomial system $\{p_{n}\}$ is orthogonal with
respect to a positive measure $\mu$, and that the convex hull of $\sp
\mu$ is a compact interval $[a, b]$. It is well known (see \cite{Ch})
that the zeros of every $p_{n}$ are contained in the interval $[a,
b]$. The following is a more precise result on the distribution of
zeros.
\medskip
{\smc Theorem A,} see \cite {Ne}, Th. 3, p. 50. Let the polynomial system
$\{p_{n}\}$ be orthogonal with respect to a compactly supported,
positive measure $\mu$ on $\bR$, and let the convex hull of $\sp \mu$ be $[a,
b]$.  Then the asymptotic distribution of zeros of $\{p_{n}\}$ is an
absolutely continuous measure $\nu$ which depends only on $[a, b]$.
The support of $\nu$ is precisely $[a, b]$ and its density in this
interval is given by
    $$\rho(x)=\frac{1}{\pi\sqrt{(b-x)(x-a)}}.$$
   \medskip

Next we describe the asymptotic distribution of zeros for a system
of eigenpolynomials.

   \medskip
   {\smc Theorem B,} see  \cite {BR} Th. 2 and 4. Let $\{p_{n}\}$ be a
   system of eigenpolynomials for an operator $\frak d$ with
   $a_{N}$ monic of degree $N$. Then the asymptotic distribution
   of zeros of $\{p_{n}\}$
   is a probability measure $\nu$ with the following properties:

   a) $\nu$ has compact support;

   b) its Cauchy transform $C(x)=\int\frac
   {d\nu(\zeta)}{x-\zeta}$ satisfies the equation
   $C(x)^{N} = 1/a_{N}(x)$ for almost all $x\in \bC$.

   These properties determine $\nu$ uniquely.
   \medskip
Note that the limiting measure $\nu$ is
independent of all terms in (2) except the leading term 
$a_{N}(x)\frac {d^{N}}{dx^N}$.
\medskip
To derive from these two results a statement about compact  
type BKS, we will need the following.
\medskip
{\smc Proposition 1.} Let $\{p_{n}\}$ be a system of eigenpolynomials
for a differential operator $\frak d$, and assume that all the zeros
of $p_{n}$ are real. Then there exists a compact set
containing all the zeros of every $p_{n}$ if and only if $\deg a_{N} =
N$.
\medskip
Now it is easy to derive the following 

\medskip
{\smc Main Theorem.} Let $\{p_{n}\}$ be a compact type BKS, orthogonal with
respect to a measure $\mu$ and with differential operator $\frak d$.
If the convex hull of $\sp \mu$ is the interval $[a, b]$, then $N$ is
even and $a_{N}(x)$ is a constant multiple of $((x-a)(x-b))^{N/2}$.

\bigskip
{\smc References and acknowledgements.} There exists a really vast
 literature devoted to the
   classification problem for OPS.
   Classification of BKS has also attracted  
   considerable attention, see
   e.g. \cite {EKLW} with its 100 references and \cite {KLY1},\cite
   {KLY2}, \cite{KYY} and references therein. The authors are happy to
   be able to contribute to this both classical and active area.
We are grateful to Dr. M.~Shapiro and Professor H.~Shapiro
for a number of discussions on the topic. The third author
wants to acknowledge the hospitality of Mathematisches
Forschungsinstitut Oberwolfach in September 2001 in whose peaceful
and serene atmosphere he  found some information on BKS and realized that the
results of \cite {MS} and \cite {BR} are
applicable to the BKS-classification problem.
\bigskip

\heading \S 2. Proofs \endheading

We need to prove Proposition 1 and the main theorem (as its corollary). 
Since the proof of Theorem B, in the situation where we will need it, 
follows easily along the same lines, we will include such a proof for 
the convenience of the reader.
\medskip
Consider a polynomial system $\{p_{n}\}$ with the associated root
measures $\nu_{n}$. Assume that the supports of the measures
$\nu_{n}$ are all contained in the same compact set, and that
$\nu_{n}  \to \nu$ in the weak topology. Let $C_{n}(x)$ be the Cauchy
transform of $\nu_{n}$ and note that
$$C_{n}(x) = \frac{p_{n}'(x)}{np_{n}(x)}.$$
If $C(x)$ denotes the Cauchy transform of $\nu$, it follows that
$$\frac{p_{n}'(x)}{np_{n}(x)} \to C(x)$$
for almost every $x \in \bC$.
\medskip
Suppose now that $\{p_{n}\}$ is a system of eigenpolynomials for an
operator $\frak d$ and that $\sp \nu_{n}$ are all contained in the
same compact subset of the real line. Then there exists at least a
subsequence of the $\nu_{n}$ converging weakly to some measure $\nu$.
Moreover, if we let $\nu_{n}^{(j)}$ denote the root measure of the
$j$th derivative of $p_{n}$, then it follows from Rolle's theorem
that (a subsequence of) $\nu_{n}^{(j)}$ converges weakly to $\nu$ for
every $j > 0$. In particular,
$$\frac{p_{n}^{(j+1)}(x)}{(n-j)p_{n}^{(j)}(x)} \to C(x)$$
for almost every $x \in \bC$, where $C(x)$ is the Cauchy transform
of $\nu$. If we divide both sides of the differential equation $\frak
d p_{n} = \lambda_{n} p_{n}$ by $n(n-1) \dots (n-N+1)p_{n}$ we obtain
$$a_{N}(x) \prod_{j =
0}^{N-1}\frac{p_{n}^{(j+1)}(x)}{(n-j)p_{n}^{(j)}(x)} +
\frac{a_{N-1}(x)}{n-N+1}\prod_{j=0}^{N-2}\frac{p_{n}^{(j+1)}(x)}{(n-j) 
p_{n}^{(j)}(x)}
+ \dots = \frac{\lambda_{n}}{n(n-1) \dots (n-N+1)}.$$
When $n \to \infty$ all the terms in the left-hand side but the first
tend to zero, and so
$$a_{N}(x)C(x)^{N} = \lim_{n \to \infty} \frac{\lambda_{n}}{n(n-1)
\dots (n-N+1)}$$ for almost all $x$.
But it can be seen (see \cite{BR}) that $\lambda_{n} =
\sum_{k=1}^{N}c_{k}n(n-1)\dots (n-k+1)$ where $c_{k}$ is the
coefficient at $x^k$ in $a_{k}(x)$. In particular, if $\deg a_{N} < 
N$ i.e. $c_{N} = 0$, then $\lambda_{n}/n(n-1) \dots (n-N+1) \to 0$
when $n \to \infty$, and it follows that $C(x) = 0$ for almost all
$x$. This implies that $\nu = 0$, a contradiction. This argument proves one of
the implications in Proposition 1. Moreover, if $a_{N}(x)$ is monic of 
degree $N$, then $c_{N}=1$ so it follows that $C(x)^{N}=1/a_{N}(x)$.  
For the other implication we refer to \cite{BR}, Lemma 9. 

\qed 
\medskip
Suppose now that $\{p_{n}(x)\}$ is a compact type BKS as in the Main 
Theorem. 
By the remark preceeding Theorem A, the zeros of every $p_{n}(x)$ are
contained in the interval $[a, b]$. It follows from Proposition 1
that $\deg a_{N}(x) = N$, and we might as well assume that $a_{N}(x)$ is
monic. Hence the Cauchy transform $C(x)$ of the limit $\nu = \lim_{n
\to \infty} \nu_{n}$ satisfies $C(x)^{N} = 1/a_{N}(x)$. On the other
hand, a direct computation of the Cauchy transform, using the
expression for $\nu$ in Theorem A, gives
$C(x)^2 = 1/(x -a)(x-b)$. Comparing these results yields
$a_{N}(x) = ((x-a)(x-b))^{N/2}$.

\qed

\bigskip
       \heading \S 3. Final remarks  \endheading

       {\smc Problem 1.} The major problem in the context of this 
       paper is whether every compact type BKS is a Jacobi-type OPS, 
       compare \cite {EKLW}, Conjecture 7.3. For the constant leading 
       coefficient the analogous fact was proved in \cite {KYY}.

      {\smc Problem 2.} Is there an analog of Theorem A on the asymptotic 
       zero distribution for a signed measure $\mu$ with compact 
       support? What is the situation for a probability measure with a noncompact 
       support as well as for Sobolev orthogonal polynomial systems. 
       (There exists a literature on this topic.) 

       {\smc Problem 3.} Generalize the results of \cite {BR} to
       operators with $\deg a_{N} < N$. Preliminary computer experiments show
       the existence of limiting measures for all operators (2).

\enddocument